\documentclass{article}
\usepackage{amssymb,latexsym,amsmath,amsthm,mathrsfs}
\usepackage{graphicx}
\newcommand{\comment}[1]{}
\begin{document}
\title{Observations on a certain theorem of Fermat and on others concerning prime numbers\footnote{Presented
to the St. Petersburg Academy on September 26, 1732.
Originally
published as
{\em Observationes de theoremate quodam Fermatiano aliisque ad numeros primos spectantibus},
Commentarii academiae scientiarum Petropolitanae \textbf{6} (1738), 103--107.
E26 in the Enestr{\"o}m index.
Translated from the Latin by Jordan Bell,
Department of Mathematics, University of Toronto, Toronto, Ontario, Canada.
Email: jordan.bell@gmail.com}}
\author{Leonhard Euler}
\date{}
\maketitle

It is known that the quantity $a^n+1$ always has divisors whenever $n$ is
an odd number or is divisible by an odd number aside from unity.\footnote{Translator: If $n=kl$ then $a^n+1$ is divisible by
\[
\frac{a^{kl}}{a^l}=a^{(k-1)l}-a^{(k-2)l}+\ldots+1.
\]
See Chapter II of Hardy and Wright.} Namely $a^{2m+1}+1$ can be divided
by $a+1$ and $a^{p(2m+1)}+1$ by $a^p+1$, for whatever number is substituted
in place of $a$. But on the other hand, if $n$ is a number which is divisible
by no odd number aside from unity, which happens when $n$ is a power
of two, 
no divisor of the number $a^n+1$ can be assigned.\footnote{Translator: Euler probably means that there is no general form for a divisor.} So if there
are prime number of this form $a^n+1$, they must all necessarily
be included in the form $a^{2^m}+1$. But it cannot however be
concluded from this that $a^{2^m}+1$ always exhibits a prime number for
any $a$;
for it is clear first that if $a$ is an odd number, this form will have
the divisor $2$. Then also, even if $a$ denotes an even number,
innumerable cases can still be given in which a composite number results.
For instance, the formula $a^2+1$ can be divided by $5$ whenever $a=5b \pm 3$,
and $30^2+1$ can be divided by $17$, and $50^2+1$ by $41$.
Similarly, $10^4+1$ has the divisor $73$, $6^8+1$ has the divisor $17$,
and $6^{128}+1$ is divisible by $257$. Yet no case has been found where
any divisor of
this form $2^{2^m}+1$ occurs, however far we have checked in the table of prime numbers, which indeed does not extend beyond $100000$.
For this and perhaps other reasons, Fermat was led to state
there to be no doubt that $2^{2^m}+1$ is always a prime number,
and proposed this eminent theorem to Wallis and other English
Mathematicians for demonstration. Indeed he admits to not himself
have a demonstration of this, but did not however hold it to be any
less than completely
true. He also praised the great utility of this, by means of which one can easily
exhibit a prime number larger than any given number, which without
a
universal theorem of this type would be very difficult. This is assembled
in the penultimate letter in the {\em Commercium Epistolicum}, included in
the
second volume of the {\em Opera} of Wallis.\footnote{Translator: See Chapter III, \S IV of Weil, {\em Number theory: an approach through history from
Hammurapi to Legendre}.} They are also recorded on
p. 115 of the works of Fermat, as follows: ``For I have said that
numbers made by 
squaring two and adding unity always lead to prime numbers, namely
that $3,5,17,257,65537$ etc. to infinity are prime, 
and the truth of this theorem has already been shown by Analysts
with no difficulty etc.''

 The truth of this theorem can be seen, as I have already said, if one
takes $1,2,3$ and $4$ for $m$; for these yields the numbers $5,7,257$ and
$65537$, which all occur among the prime numbers in the table. But I do not know
by what fate it turned out that the number immediately following,
$2^{2^5}+1$, ceases to be a prime number; for I have observed after thinking about
this for many days that this number can be divided by $641$,
which can be seen at once by anyone who cares to check.
For it is $2^{2^5}+1=2^{32}+1=4294967297$.
From this it can be understood that the theorem fails in this and even in other
cases which follow, and hence the problem of finding a prime number greater
than a given number still remains unsolved.

I will now examine also the formula $2^n-1$, which, whenever $n$ is not
a prime number, has divisors, and not only $2^n-1$, but also $a^n-1$. But if
$n$ is a prime number, it might seem that $2^n-1$ also always gives a prime;
this however no one, as far as I know, has dared to profess, and indeed it can
easily be refuted. Namely $2^{11}-1$, i.e. $2047$, has the divisors
$23$ and $89$, and
$2^{23}-1$ can be divided by $47$. I see also that the Cel. Wolff has not
only not mentioned this in the new edition of his {\em Elem. Matheseos},
where he investigates the perfect numbers and includes $2047$ among the primes,
but also has $511$ or $2^9-1$ as a prime, while it is divisible by $2^3-1$, i.e.
$7$. He also gives that $2^{n-1}(2^n-1)$ is a perfect number whenever
$2^n-1$ is prime; therefore $n$ must also be a prime number. 
I have found it a worthwhile effort to examine those cases in which $2^n-1$
is not a prime number while $n$ is. I have also found that
$n=4m-1$ and $8m-1$ are prime numbers, then $2^n-1$ can always be divided by
$8m-1$. Hence the following cases should be excluded: $11$, $23$, $83$, $131$,
$179$, $191$, $239$ etc., which numbers when substituted for $n$
yield $2^n-1$ that is a composite number.
Neither however can all the remaining prime numbers be successfully
put in place of $n$, but 
still more must be removed; thus I have observed that $2^{37}-1$ can be
divided by $223$, $2^{43}-1$ by $431$,
$2^{29}-1$ by $1103$, $2^{73}-1$ by $439$;
however it is not in our power to exclude them all. Still,
I venture to assert that except for those cases noted, all prime numbers less than $50$ and perhaps even $100$ yield $2^{n-1}(2^n-1)$ which is a perfect number,
thus $11$ perfect numbers arise from the following numbers taken for $n$,
$1$, $2$, $3$, $5$, $7$, $13$, $17$, $19$, $31$, $41$, $47$.
I have deduced these observations from a not inelegant theorem, whose
proof I do not have, but indeed of whose truth I am completely certain.
This theorem is: {\em $a^n-b^n$ can always be divided by $n+1$, if $n+1$
is any prime number which divides neither $a$ or $b$}; I believe this
demonstration is more difficult because it is not true unless $n+1$ is
a prime number. From this theorem, it follows at once that $2^n-1$ can
always be divided by $n+1$ if $n+1$ is a prime number, or, 
since each prime aside from $2$ is odd, and as $a=2$ that case does not happen
because of the conditions of the theorem,\footnote{Translator: If $n+1=2$ then
$n+1$ divides $a$.} $2^{2m}-1$ will always be able to be divided by $2m+1$
if $2m+1$ is a prime number.
Hence either $2^m+1$ or $2^m-1$ will
be able to be divided by $2m+1$.\footnote{Translator: $2^{2m-1}=(2^m+1)(2^m-1)$,
and since $2m+1$ is prime and divides $2^{2m-1}$, it must divide one of the factors.} I have also discovered that $2^m+1$ can be divided if
$m=4p+1$ or $4p+2$; while $2^m-1$ will have the divisor $2m+1$ if $m=4p$ 
or $4p-1$. I have happened upon many other theorems in this pursuit which
are no less elegant,
which I believe should be further investigated, because either they
cannot be demonstrated themselves, or they follow from propositions which cannot be
demonstrated; some which seem important are appended here.

\begin{center}
{\Large Theorem 1}
\end{center}

{\em If $n$ is a prime number, all powers having the exponent $n-1$ leave
either nothing or $1$ when divided by $n$.}

\begin{center}
{\Large Theorem 2}
\end{center}

{\em With $n$ still a prime number, every power whose exponent is $n^{m-1}(n-1)$
leaves either $0$ or $1$ when divided by $n^m$.}

\begin{center}
{\Large Theorem 3}
\end{center}

{\em Let $m,n,p,q$ etc be distinct prime numbers and let $A$ be the least
common multiple of them decreased by unity, think of them $m-1,n-1,p-1,q-1$ etc.;
with this done, I say that any power of the exponent $A$, like $a^A$,
divided by $mnpq$ etc. will leave either $0$ or $1$, unless
$a$ can be divided by one of the numbers $m,n,p,q$ etc.}

\begin{center}
{\Large Theorem 4}
\end{center}

{\em With $2n+1$ denoting a prime number, $3^n+1$ will be able to be divided
by $2n+1$, if either $n=6p+2$ or $n=6p+3$;
while $3^n-1$ will be able to be divided by $2n+1$ if either
$n=6p$ or $n=6p-1$.}

\begin{center}
{\Large Theorem 5}
\end{center}

{\em $3^n+2^n$ can be divided by $2n+1$ if $n=12p+3$, $12p+5$, $12p+6$
or $12p+8$, And $3^n-2^n$ can be divided by $2n+1$ if $n=12$, $12p+2$,
$12p+9$ or $12p+11$.}

\begin{center}
{\Large Theorem 6}
\end{center}

{\em Under the same conditions which held for $3^n+2^n$, $6^n+1$ can also
be divided by $2n+1$; and $6^n-1$ under those which held for $3^n-2^n$.}

\end{document}